\documentclass[12pt]{article}
\usepackage{amsmath}
\usepackage{amssymb}
\usepackage{amsthm}
\usepackage{amscd}
\usepackage{color}
\usepackage[mathscr]{eucal}
\usepackage[all]{xy}

\DeclareMathOperator{\Hom}{Hom}

\DeclareMathOperator{\End}{End}

\DeclareMathOperator{\id}{id}

\DeclareMathOperator{\Real}{Re}

\newcommand{\C}{\mathbb{C}}
\newcommand{\Z}{\mathbb{Z}}

\newcommand{\N}{\mathbb{N}}

\newcommand{\CP}{\mathbb{P}}

\newcommand{\vac}{\mathbf{1}}

\newcommand{\frakg}{\mathfrak{g}}

\newcommand{\NO}{\,{\raise0.25em\hbox{$\mathop{\hphantom{\cdot}}%
\limits^{_{\circ}}_{^{\circ}}$}}\,}
\newcommand{\fusion}[3]{{\binom{#3}{#1\;#2}}}

\theoremstyle{plain}
   \newtheorem{theorem}
{Theorem}
[subsection]
   \newtheorem{corollary}[theorem]{Corollary}
   \newtheorem{lemma}[theorem]{Lemma}
   \newtheorem{proposition}[theorem]{Proposition}
\theoremstyle{remark}
   \newtheorem{definition}[theorem]{Definition}

\numberwithin{equation}{section}

\begin{document}
\begin{center}
\begin{large}
Logarithmic intertwining operators \\
and\\
the space of conformal blocks over
the projective line
\end{large}
\end{center}
\vskip5ex
\begin{center}
Yusuke Arike
\vskip1ex
Department of Mathematics,
Graduate School of Science\\
Osaka University\\
y-arike@cr.math.sci.osaka-u.ac.jp
\end{center}
\vskip2ex
\abstract{
We show that the space of logarithmic intertwining operators among logarithmic modules
for a vertex operator algebra is isomorphic to 
the space of $3$-point conformal blocks over the projective line.
This is considered as a generalization of Zhu's result for ordinary intertwining operators
among ordinary modules.}
\vskip2ex

\section{Introduction}
One of the most important problems in representation theory of vertex operator algebras
is to determine fusion rules which are the dimensions of intertwining operators among
three modules for vertex operator algebras.
Intertwining operators of the type $\fusion{M^1}{M^2}{M^3}$
are linear maps $I(-,z):M^1\to\Hom_\C(M^2,M^3)[[z,z^{-1}]]$ with several axioms (see \cite{FHL})
where $M^i\,(i=1,2,3)$ are modules for a vertex operator algebra.

The definition of intertwining operators given in \cite{FHL} treats modules on which $L_0$ acts as a semisimple operator.
However, in general, we have to consider
modules which do not decompose into $L_0$-eigenspaces
but do into generalized $L_0$-eigenspaces
Such modules are called {\it logarithmic modules} in \cite{M1}.

A notion of {\it logarithmic intertwining operators} among logarithmic modules is introduced in \cite{M1}.
Logarithmic intertwining operators may involve logarithmic terms.
It is shown in \cite{M1} that a logarithmic intertwining operator among
ordinary modules is nothing but the so-called intertwining operator.
Several examples of logarithmic modules are found and
logarithmic intertwining operators among these modules are constructed (see eg. \cite{M1}, \cite{M2}, \cite{AM}).

On the other hand in conformal field theory its important feature is a notion 
of conformal blocks associated with vertex operator algebras.
Mathematically rigorous formulation of $N$-point conformal blocks on Riemann
surfaces associated with vertex operator algebras
is given in \cite{Z1} with the assumption that the corresponding vertex operator algebra is quasi-primary generated.
It is shown in \cite{Z1} that the space of $3$-point conformal blocks over the projective
line $\CP^1$ is isomorphic to the space of intertwining operators among ordinary modules
for a vertex operator algebra.

In this paper we give a sort of generalization of Zhu's result in the case that
the modules are logarithmic.
More precisely we are going to prove that the space of $3$-point conformal blocks over the projective line
is isomorphic to the space of logarithmic intertwining operators without the assumption that
a vertex operator algebra is quasi-primary generated.
Taking the formulation of the space of coinvariants in \cite{NT} we do not have to assume that a vertex operator algebra is quasi-primary generated.

The study on logarithmic intertwining operators is very important since if we could know its dimension from $S$-matrix obtained by formal characters
in fact if a vertex operator algebra is rational and satisfies several conditions
the dimension of intertwining operators is completely determined by $S$-matrix.
However it is left for further studies.

The paper is organized as follows.
In section 2 we recall the definition of vertex operator algebras and their modules.
The definition of logarithmic modules is located here.
Also we describe the space of conformal blocks over $\CP^1$ according to \cite{NT}.

In section 3 we recall the definition and properties of logarithmic intertwining operators
which are given in \cite{M1} and \cite{HLZ}.

We state the main theorem in this paper and give a proof in section 4.
The linear maps between the space of logarithmic intertwining operators
and the space of $3$-point conformal blocks are defined and it is proved that these maps are 
well-defined and mutually inverse.

\section{Vertex operator algebras and the space of conformal blocks over the projective line}
Throughout this paper we use the notation $\N=\{0,1,2,\dots\}$.
\subsection{Vertex operator algebras and current Lie algebras}
A {\it vertex operator algebra} is a $\N$-graded vector space $V=\bigoplus_{k=0}^\infty V_k$
with $\dim V_k<\infty \,(k\in\Z_{\ge0})$ equipped with a linear map
\begin{equation}
Y(-,z):V\to\End(V)[[z,z^{-1}]],\,Y(a,z)=\sum_{n\in\Z}a_{(n)}z^{-n-1}
\end{equation}
and with distinguished vectors $\vac\in V_0$ called the {\it vacuum vector}
and $\omega\in V_2$ called the {\it Virasoro vector} satisfying the following axioms (see e.g. \cite{FHL}, \cite{MN}):
\vskip 1ex
\noindent
(1) For any pair of vectors in $V$ there exists a nonnegative integer $N$ such that
$a_{(n)}b=0$ for all integers $n\ge N$.\\
\noindent
(2)
For any vectors $a, b, c\in V$ and integers $p, q, r\in\Z$, 
\begin{equation}
\begin{split}
&\sum_{i=0}^\infty\binom{p}{i}(a_{(r+i)}b)_{(p+q-i)}c\\
&\qquad\qquad=\sum_{i=0}^\infty(-1)^i\binom{r}{i}(a_{(p+r-i)}b_{(q+i)}c-(-1)^rb_{(q+r-i)}a_{(p+i)}c)
\end{split}
\end{equation}
hold.\\
\noindent
(3)
$Y(\vac,z)=\id_V$.\\
(4)
$Y(a,z)\vac\in V[[z]]$ and $a_{(-1)}\vac=a$.\\
\noindent
(5)
Set $L_n=\omega_{(n+1)}$.
Then $\{L_n\,|\,n\in\Z\}$ together with the identity map on $V$ give a representation of the Virasoro algebra on $V$ with central charge $c_V\in\C$.\\
\noindent
(6)
$L_0a=ka$ for any $a\in V_k$ and nonnegative integers $k$.

\vskip 1ex
\noindent
(7)
$\dfrac{d}{dz}Y(a,z)=Y(L_{-1}a,z)$ for any $a\in V$.

\medskip
\noindent
(8)
Denote $|a|=k$ for any $a\in V_k$ and then
\begin{equation}
|a_{(n)}b|=|a|+|b|-1-n
\end{equation}
for any homogeneous $b\in V$ and $n\in\Z$.

\medskip

In order to define the space of conformal blocks we introduce the spaces 
$V^{(1)}=\bigoplus_{k=0}^\infty V_k\otimes \C((\xi))(d\xi)^{1-k}$
and $V^{(0)}=\bigoplus_{k=0}^\infty V_k\otimes\C((\xi))(d\xi)^{-k}$.
Let $\nabla:V^{(0)}\to V^{(1)}$ be the linear map defined by
\begin{equation}
v\otimes f(\xi)(d\xi)^{-n}=L_{-1}v\otimes f(\xi)(d\xi)^{-k}+v\otimes\frac{df(\xi)}{d\xi}(d\xi)^{1-k}.
\end{equation}
We set $\frakg=V^{(1)}/\nabla V^{(0)}$
and denote  the image of $a\otimes f(\xi)(d\xi)^{1-k}\in V_k\otimes\C((\xi))(d\xi)^{1-k}$ by $J(a,f)$. Then we have:

\begin{proposition}[\cite{NT}, Proposition 2.1.1]\label{prop-current}
The vector space $\frakg$ is a Lie algebra with the blacket
\begin{equation}
[J(a,f), J(b,g)]=\sum_{m=0}^{|a|+|b|-1}\frac{1}{m!}J\bigl(a_{(m)}b, \frac{d^mf}{d\xi^m}g\bigr)
\end{equation}
for homogeneous $a, b\in V$.
\end{proposition}

The Lie algebra $\frakg$ is called the {\it current Lie algebra}.
Let us denote $J_n(a)=J(a,\xi^{n+|a|-1})$.

Applying the construction of the current Lie algebra
to the vector space $\bigoplus_{k=0}^\infty V_k\otimes\C[\xi,\xi^{-1}](d\xi)^{1-k}$,
we have a graded Lie algebra $\bar{\frakg}=\bigoplus_{n\in\Z}\bar{\frakg}_n$
where the vector space $\bar{\frakg}_n$ is linearly spanned by $J_n(a)\,(a\in V)$. 
The Lie algebra $\bar{\frakg}$ is a Lie subalgebra of $\frakg$.
The following proposition plays an important role when we define duality functor
on the category of $V$-modules.

\begin{proposition}[\cite{NT}, Proposition 4.1.1]\label{prop-dual}
The linear map $\theta:\bar{\frakg}\to\bar{\frakg}$ defined by
\begin{equation}\label{eq-involution}
\theta(J_n(a))=(-1)^{|a|}\sum_{j=0}^\infty\frac{1}{j!}J_{-n}(L_1^ja)
\end{equation}
for $a\in V$ and $n\in\Z$ is an anti-Lie algebra involution.
\end{proposition}

\subsection{Modules for vertex operator algebras}
Let $M$ be a weak $V$-module (see \cite{DLM} for the definition).
A weak $V$-module $M$ is called $\N$-gradable if it admits a decomposition
$M=\bigoplus_{n\in\N}M_{(n)}$ such that
\begin{equation}
a_{(n)}M_{(k)}\subseteq M_{(|a|+k-1-n)}
\end{equation}
for homogeneous $a\in V$ and $n\in\Z$.

Let $M$ be a weak $V$-module.
A weak $V$-module is called a {\it logarithmic module}
if $M$ decomposes into a direct sum of generalized $L_0$-eigenspaces.

Let $M=\bigoplus_{n\in\Z_{\ge0}}M_{(h+n)}$ be a logarithmic module with a complex number $h$ and
\begin{equation}
M_{(h+n)}=\{x\in M\,|\,(L_0-h-n)^{k+1}x=0 \text{ for a nonnegative integer }k
\}.
\end{equation}
Obviously $M$ is a $\N$-gradable $V$-module.

In this paper a $V$-module $M$ is always a logarithmic $V$-module 
satisfies the following conditions.

\medskip
\noindent
i) There exist a complex number $h$ and a nonnegative integer $k$ such that
$M=\bigoplus_{n=0}^\infty M_{(h+n)}$ with
$M_{(h+n)}=\{u\in M\,|\,(L_0-h-n)^{k+1}u=0\}$ for all $n\in\Z$.\\
\noindent
ii) $\dim M_{(h+n)}<\infty$ for all nonnegative integers $n$.
We denote $|u|=h+n$ for $u\in M_{(h+n)}$ for short.

\medskip
We remark that any $V$-module in this paper is a $V$-module in the sense of  \cite{NT} and \cite{MNT}.

Let $k$ be a nonnegative integer and let $\mathscr{C}_k$ be the category consisting of $V$-modules whose homogeneous subspaces
are annihilated by $(L_0-h-n)^{k+1}$.
Then it follows that $\mathscr{C}_0\subseteq \mathscr{C}_1\subseteq\dotsb
\subseteq\mathscr{C}_k\subseteq\dotsb$.

Any $V$-module $M$ is a $\bar{\frakg}$-module by the action
\begin{equation}\label{eq-action}
J_n(a)u=a_{(|a|-1+n)}u 
\end{equation}
for any homogeneous $a\in V$ and $u\in M$ (cf. \cite{DLM}, \cite{NT}).
For any $a\in V$ and $u\in M$, there exists a nonnegative integer $n_0$ such that
$a_{(n)}u=0$ for all $n\ge n_0$. Therefore, the $V$-module
$M$ is also a $\frakg$-module by the action \eqref{eq-action}.

Let us denote the restricted dual of a $V$-module $M=\bigoplus_{n=0}^\infty M_{(h+n)}$ by $D(M)=\bigoplus_{n=0}^\infty M_{(h+n)}^\ast$
where $M_{(h+n)}^\ast=\Hom_\C(M_{(h+n)},\C)$.
A $\bar{\frakg}$-module structure on $D(M)$ can be defined by letting
\begin{equation}
\langle J_n(a)\varphi, u\rangle=\langle \varphi, \theta(J_n(a))u\rangle
\end{equation}
for all $\varphi\in D(M)$ and $u\in M$.
The following proposition is known.

\begin{proposition}[\cite{NT}, Proposition 4.2.1, cf. \cite{FHL}, Theorem 5.2.1]
There exists a unique $V$-module structure on $D(M)$
which extends $\bar{\frakg}$-module structure.
\end{proposition}

Since $\langle L_n\varphi,u\rangle=\langle \varphi, L_{-n}u\rangle$ for all $\varphi\in D(M)$ and $u\in M$,
we see that $D(M)=\oplus_{n=0}^{\infty}D(M)_{(h+n)}$ and $D(M)
\in\mathscr{C}_k$
for any $M\in\mathscr{C}_k$.

\subsection{The space of conformal blocks over the projective line}

Let $\CP^1=\C\cup\{\infty\}$ be the projective line
and $z$ its inhomogeneous coordinate.
Let $A=\{1,2,\dotsc,N,\infty\}$ and let us fix a set $p_A=(p_a)_{a\in A}$ of $N+1$ distinct points $p_a\in\CP^1\,(a\in A)$
with $p_\infty=\infty$.
We write $\xi_a=z-p_a\,(a\not=\infty)$ and $\xi_\infty=z$, respectively.

We denote by $H^0(\CP^1,\Omega^{1-k}(\ast p_A))$
the vector space of global meromorphic $(1-k)$-differentials whose poles are located only at $p_a\,(a\in A)$.
Set $H(V,\ast p_A)^{(1)}=\bigoplus_{k=0}^\infty V_k\otimes H^0(\CP^1,\Omega^{1-k}(\ast p_A))$
and $H(V,\ast p_A)^{(0)}=\bigoplus_{k=0}^\infty V_k\otimes H^0(\CP^1,\Omega^{-k}(\ast p_A))$.
Define the linear map $\nabla:H(V,\ast p_A)^{(0)}\to H(V,\ast p_A)^{(1)}$ by
\begin{equation}
a\otimes f(z)(dz)^{1-k}\mapsto L_{-1}a\otimes f(z)(dz)^{-k}+a\otimes \frac{df(z)}{dz}(dz)^{1-k}\qquad(a\in V_k).
\end{equation}
We set 
\begin{equation}
\frakg(\CP^1,\ast p_A)=H(V,\ast p_A)^{(1)}/\nabla H(V,\ast p_A)^{(0)}.
\end{equation}
It is shown (cf. \cite[Proposition 5.1.1]{NT}) that the vector space 
$\frakg(\CP^1,\ast p_A)$
is a Lie algebra with the blacket
\begin{multline}
\bigl[a\otimes f(z)(dz)^{1-|a|}, b\otimes g(z)(dz)^{1-|b|}\bigr]\\
=\sum_{m=0}^\infty\frac{1}{m!}a_{(m)}b\otimes \frac{d^mf(z)}{dz^m}g(z)(dz)^{2-|a|-|b|+m}.
\end{multline}

For each $a\in A$ we define the linear map
\begin{equation}
i_a:H^0(\CP^1,\Omega^k(\ast p_A))\to
\begin{cases}
\C((\xi_a))(d\xi_a)^k, & a\in A\backslash\{\infty\}\\
&\\
\C((\xi_\infty^{-1}))(d\xi_\infty)^k, & a=\infty
\end{cases}
\end{equation}
by taking the Laurent expansion at $z=p_a$ in terms of the coordinate $\xi_a$.
We denote $i_a f(z)(dz)^k$ by $f_a(\xi_a)(d\xi_a)^k$.

For any $a\in A\backslash\{\infty\}$, we define the linear map 
$j_a:\frakg(\CP^1,\ast p_A)\to\frakg$
by $j_a(a\otimes f(z)(dz)^{1-k})=a\otimes f_a(\xi_a)(d\xi_a)^{1-k}$ 
and the linear map $j_\infty:\frakg(\CP^1,\ast p_A)\to\frakg$
by $j_\infty(a\otimes f(z)(dz)^{1-k})=-\theta(a\otimes f_\infty(\xi_\infty)(d\xi_\infty)^{1-k})$.
Then the linear map $j_\infty$ is well-defined
since
\begin{equation}
j_\infty(a\otimes f(z)(dz)^{1-k})=-\sum_{n\le n_0}f_n\theta(J_n(a)),\;
\theta(J_n(a))=(-1)^kJ_{-n}(e^{L_1}a),
\end{equation}
where $f_\infty(\xi_\infty)=\sum_{n\le n_0}f_n \xi_\infty^{n+k-1}$ (see \cite{NT}).
The following proposition is fundamental.

\begin{proposition}[\cite{NT}, Proposition 5.1.3]\label{prop-Lie}
For any $a\in A$, the linear map $j_a:\frakg(\CP^1,\ast p_A)\to\frakg$
is a Lie algebra homomorphism.
\end{proposition}

Let $M^a\,(a\in A)$ be $V$-modules.
We set $M_A=\bigotimes_{a\in A}M^a$ and $\frakg_A=\frakg^{\oplus|A|}$.
Let $\rho_a:\frakg\to\End(M^a)$ be the representation defined by \eqref{eq-action}
for $a\in A$.
Then the linear map $\rho_A:\frakg_A\to\End(M_A)$ defined by $\rho_A=\oplus_{a\in A}\rho_a$
is a representation of the Lie algebra $\frakg_A$ on $M_A$.
We denote the image of the Lie algebra homomorphism $j_A=\sum_{a\in A}j_a$
by $\frakg_{p_A}^{out}$, which acts on $M_A$ via $\rho_A$.
The following definition is given by \cite{NT}.

\begin{definition}
The vector space $C^*(M_A, p_A)=\Hom_\C(M_A/\frakg_{p_A}^{out}M_A,\C)$
is called the space of conformal blocks at $p_A$.
\end{definition}

\section{Logarithmic intertwining operators}
In this section, we recall the notion of logarithmic intertwining operators 
and their properties according to \cite{M1}.
\subsection{Definition}

\begin{definition}[\cite{M1}]
Let $M^1,\,M^2$ and $M^3$ be weak $V$-modules.
A {\it logarithmic intertwining operator of the type $\fusion{M^1}{M^2}{M^3}$} is a linear map
\begin{align}
&I(-,z):M^1\to\Hom_\C(M^2,M^3)\{z\}[\log z]\\
&I(u,z)=\sum_{n=0}^d\sum_{\alpha\in\C}u_{(\alpha)}^nz^{-\alpha-1}(\log z)^k
\end{align}
with the following properties:\\
\noindent
i)\,(Truncation condition)
For any $u_1\in M^1$, $u_2\in M^2$ and $0\le k\le d$,
\begin{equation}
(u_1)_{(\alpha)}^ku_2=0
\end{equation}
for sufficiently large $\Real(\alpha)$.\\
\noindent
ii)\,($L_{-1}$-derivative property)
For any $u_1\in M^1$,
\begin{equation}
I(L_{-1}u_1,z)=\frac{d}{dz}I(u_1,z).
\end{equation}
\noindent
iii)\,
For all $a\in V$, $u_1\in M_1$, $u_2\in M_2$, $\alpha\in\C$, $0\le n\le d$ and $p,q\in\Z$, we have
\begin{equation}\label{eq-borcherds}
\begin{split}
&\sum_{i=0}^\infty\binom{p}{i} (a_{(q+i)}u_1)_{(\alpha+p-i)}^n\\
&\qquad\qquad=\sum_{i=0}^\infty(-1)^i\binom{q}{i}(a_{(p+q-i)}(u_1)^n_{(\alpha+i)}-(-1)^q(u_1)_{(\alpha+q-i)}^na_{(p+i)}). 
\end{split}
\end{equation}

We denote the space of logarithmic intertwining operators of the type $\fusion{M^1}{M^2}{M^3}$
by $I\fusion{M^1}{M^2}{M^3}$, that is, we use the same notation as usual intertwining operators.
\end{definition}

\medskip
Setting $q=0$ and $p=0$ in \eqref{eq-borcherds}, respectively, we have
\begin{align}
&[a_{(p)},(u_1)_{(\alpha)}^n]=\sum_{i=0}^\infty \binom{p}{i}(a_{(i)}u_1)_{(\alpha+p-i)}^n,\label{eq-com}\\
&(a_{(q)}u_1)_{(\alpha)}^n=\sum_{i=0}^\infty(-1)^i\binom{q}{i}\label{eq-ass}
\{a_{(q-i)}(u_1)_{(\alpha+i)}^n-(-1)^qa_{(\alpha+q-i)}^na_{(i)}\}
\end{align}
and we call, by abuse of terminologies,  the {\it commutator formula} and {\it associativity formula}, respectively.
By the commutator formula, we have
\begin{equation}
[L_{-1}, u_{(\alpha)}^n]=(L_{-1}u)_{(\alpha)}^n\label{eq-derivation}
\end{equation}
for any $u\in M^1$ and $0\le n\le d$.
By the associativity formula, \eqref{eq-derivation} and $L_{-1}$-derivative property,
we have
\begin{equation}
(L_0u)_{(\alpha)}^n=
\begin{cases}
[L_0, (u)_{(\alpha)}^n]+(\alpha+1)(u)_{(\alpha)}^n-(n+1)
(u)_{(\alpha)}^{n+1} & 0\le n\le d-1,\\
&\\
[L_0, (u)_{(\alpha)}^n]+(\alpha+1)(u)_{(\alpha)}^n & n=d\\
\end{cases}
\label{eq-fund}
\end{equation}
for any $u\in M^1$.

\subsection{Properties for logarithmic intertwining operators}

Let $M^i=\bigoplus_{n=0}^\infty M^i_{(h_i+n)}\,(i=1,2,3)$ be objects in $\mathscr{C}_{k_i}$
for nonnegative integers $k_i \,(i=1,2,3)$ 
and complex numbers  $h_i\,(i=1,2,3)$.
Suppose that a logarithmic intertwining operator $I(-,z)$ of the type $\fusion{M^1}{M^2}{M^3}$ is
of the form
\begin{equation}
I(u_1,z)=\sum_{n=0}^d\sum_{\alpha\in\C}(u_1)_{(\alpha)}^nz^{-\alpha-1}(\log z)^n.
\end{equation}
For any homogeneous element $u_i\in M^i\,(i=1,2)$, we introduce notations
\begin{align}
x_1(u_1)_{(\alpha)}^nu_2&=((L_0-|u_1|)u_1)_{(\alpha)}^nu_2,\\
x_2(u_1)_{(\alpha)}^nu_2&=(u_1)_{(\alpha)}^n(L_0-|u_2|)u_2,\\
x_3(u_1)_{(\alpha)}^nu_2&=(L_0+\alpha+1-|u_1|-|u_2|)(u_1)_{(\alpha)}^nu_2.
\end{align}
Note that these operations $x_1$ and $x_2$ are mutually commutative (see \cite{M1}).
By using these operations we get:

\begin{lemma}[\cite{HLZ}, Lemma 3.8]\label{lemma-fund}
Let 
\begin{equation}
I(-,z)=\sum_{n=0}^d\sum_{\alpha\in\C}(-)_{(\alpha)}^nz^{-\alpha-1}(\log z)^n
\end{equation}
be a logarithmic intertwining operator of the type $\fusion{M^1}{M^2}{M^3}$
and let $p, q$ be integers such that $p\ge0$ and $0\le q\le d$.
Then
\[
x_3^p(u_1)_{(\alpha)}^qu_2=\sum_{\ell=0}^N\binom{p}{\ell}\frac{(q+\ell)!}{q!}
(x_1+x_2)^{p-\ell}(u_1)_{(\alpha)}^{q+\ell}u_2
\]
for homogeneous $u_1\in M^1$ and $u_2\in M^2$ where 
$N=\min\{p,\,d-q\}$.
\end{lemma}

The following proposition is proved in \cite{M1} by using differential equations
and in \cite[Proposition 3.9]{HLZ} by using Lemma \ref{lemma-fund}.

\begin{proposition}[\cite{M1}, Proposition 1.10]\label{prop-eigen}
Suppose that $M^i\in\mathscr{C}_{k_i}\,(i=1,2,3)$ for nonnegative integers $k_i$
and that $M^i=\bigoplus_{n=0}^\infty M^i_{(h_i+n)}$ for complex numbers $h_i\,(i=1,2,3)$.
Let $I(-,z)\in I\fusion{M^1}{M^2}{M^3}$ be a logarithmic intertwining operator
such that
\begin{equation}
I(u_1,z)=\sum_{n=0}^d\sum_{\alpha\in\C}(u_1)_{(\alpha)}^n z^{-\alpha-1}(\log z)^n \qquad (u_1\in M^1).
\end{equation}
\noindent
{\rm(1)}\,For any homogeneous $u_i\in M^i\,(i=1,2)$ we have
$|(u_1)_{(\alpha)}^nu_2|=|u_1|+|u_2|-1-\alpha$
for all $0\le n\le d$.\\
\noindent
{\rm(2)}\, For any $u_i\in M^i\,(i=1,2)$ we have
\[
I(u_1,z)u_2\in\sum_{n=0}^{k_1+k_2+k_3}M^3((z))z^{h_3-h_1-h_1}(\log z)^n.
\]
\end{proposition}
%

\section{The space of $3$-point conformal blocks and logarithmic intertwining operators}
In this section we focus on $3$-point conformal blocks in conformal field theories
over the projective line.
We prove that the space of $3$-point conformal blocks over $\CP^1$ is isomorphic to the space of logarithmic intertwining operators.
The almost same result is found in \cite{Z1}, however, the categories of modules
of us and the one in \cite{Z1} are slightly different.
\subsection{Main theorem}
Set $A=\{1,2,\infty\}$ and
let $p_A=\{0,1,\infty\}$ be the set of points on $\CP^1$.
Let $z$ be the inhomogeneous  coordinate of $\CP^1$.
The $\xi_0=z$, $\xi_1=z-1$ and $\xi_\infty=z$ are local coordinate of $\CP^1$ at $0,1,$ and $\infty$, respectively.
Take $V$-modules $M^1,M^2$ and $M^3$.
We assume that there exist complex numbers $h_i\in\C \,(i=1,2,3)$ such that
$M^i=\bigoplus_{n=0}^\infty M_{(h_i+n)}$
and that $M^i\in\mathscr{C}_{k_i}\,(i=1,2,3)$ for nonnegative integers $k_i \,(i=1,2,3)$.
Let us set $M_A=M^1\otimes M^2\otimes M^3$.
We denote the space of conformal blocks at $p_A=\{0,1,\infty\}$
by $C^*(M_A,p_A)$.
Then we can now state the main theorem of the paper which is
a generalization of Zhu's result \cite[Proposition 7.4]{Z1}.

\begin{theorem}\label{thm-1}
Let $M^i (i=1,2,3)$ be $V$-modules with
$M^i=\bigoplus_{n=0}^\infty M^i_{(h_i+n)}$ and 
 $M^i\in\mathscr{C}_{k_i}\,(i=1,2,3)$ for nonnegative integers $k_i \,(i=1,2,3)$.
The space of conformal blocks $C^*(M_A,p_A)$ at $p_A=\{0,1,\infty\}$
is isomorphic to the space of logarithmic intertwining operators of the type 
$\fusion{M^2}{M^1}{D(M^3)}$
\end{theorem}

Let $C_2(V)$ be the vector subspace of $V$ spanned by vectors of the form
$a_{(2)}b\,(a,b\in V)$.
If $\dim V/C_2(V)<\infty$ we say that $V$ satisfies {\it Zhu's finiteness condition}
which is introduced in \cite{Z2}.

By combining \cite[Theorem 5.8.1]{NT} and the theorem we get:
\begin{corollary}
If $V$ satisfies Zhu's finiteness condition then the space of intertwining operators is finite-dimensional.
\end{corollary}

\subsection{Proof of Theorem \ref{thm-1}}
For any logarithmic intertwining operator $I(-,z)$ of the type $\fusion{M^2}{M^1}{D(M^3)}$,
we define $F\in\Hom_\C(M_A,\C)$ by
\begin{equation}\label{eq-main-def-1}
\left\langle F, u_1\otimes u_2\otimes u_3\right\rangle=\left\langle I(u_2,1)u_1,u_3\right\rangle
\end{equation}
for any $u_1\in M^1$, $u_2\in M^2$ and $u_3\in M^3$.
For any $V$-module $M\in\mathscr{C}_{k}$
we define the operator $z^{L_0}:M\to M\{z\}[\log z]$ by
\begin{equation}
z^{L_0}u=\sum_{j=0}^{k}\frac{1}{j!}(L_0-|u|)^jz^{|u|}(\log z)^j.
\end{equation}
For any $x\in C^*(M_A,p_A)$, we define $I_x(-,z)\in\Hom_\C(M^1, D(M^3))\{z\}[\log z]$ by
\begin{multline}\label{eq-main-def-2}
\left\langle I_x(u_2,z)u_1,u_3\right\rangle\\
=\left\langle x, z^{-L_0}u_1\otimes z^{-L_0}u_2\otimes z^{L_0}u_3\right\rangle
\text{ for all }u_i\in M^i\,(i=1,2,3).
\end{multline}

We are going to give a prove of the theorem by dividing its into three steps.
In step 1 we prove that $F$ belongs to $C^*(M_A,p_A)$ and
show that $I_x$ is a logarithmic intertwining operator
among $V$-modules in step 2.
The final step is devoted to the proof that the correspondence between $F$ and $I_x$ is one-to-one.

\medskip
\noindent
{\bf (Step 1)} 
In order to prove that $F$ belongs to $C^*(M_A,p_A)$,
by the definition of the space of conformal blocks, it is sufficient to prove
that
\begin{equation}\label{eq-main-0}
\begin{split}
&\left\langle F, j_0(a\otimes f(z)(dz)^{1-k})u_1\otimes u_2\otimes u_3\right\rangle\\
&\quad\qquad\qquad+\left\langle F, u_1\otimes j_1(a\otimes f(z)(dz)^{1-k})u_2\otimes u_3\right\rangle\\
&\qquad\qquad\qquad\qquad\qquad+\left\langle F, u_1\otimes u_2\otimes j_\infty(a\otimes f(z)(dz)^{1-k})u_3\right\rangle
=0
\end{split}
\end{equation}
for all $a\in V_k$ and $f(z)(dz)^{1-k}\in H^0(\CP^1,\Omega^{1-k}(\ast p_A))$.
It is well known that  $\{z^p(z-1)^q(dz)^{1-k}\,|\,p,q\in\Z\}$ is a topological basis of $H^0(\CP^1,\Omega^{1-k}(\ast p_A))$.
Therefore, it is enough to show \eqref{eq-main-0} for 
$f(z)=z^p(z-1)^q\,(p,q\in\Z)$.
First of all we have
\begin{equation}\label{eq-main-1}
\begin{split}
j_0(a\otimes z^p(z-1)^q(dz)^{1-k})u_1&=\Bigl(\sum_{i=0}^\infty(-1)^{q-i}\binom{q}{i}a\otimes \xi_0^{p+i}(d\xi_0)^{1-k}\Bigr)u_1\\
&=\sum_{i=0}^\infty(-1)^{q-i}\binom{q}{i}J_{p+i-k+1}(a)u_1\\
&=\sum_{i=0}^\infty(-1)^{q-i}a_{(p+i)}u_1
\end{split}
\end{equation}
and secondly
\begin{equation}\label{eq-main-2}
\begin{split}
j_1(a\otimes z^p(z-1)^q(dz)^{1-k})u_2&=\Bigl(\sum_{i=0}^\infty\binom{p}{i}a\otimes \xi_1^{q+i}(d\xi_1)^{1-k}\Bigr)u_2\\
&=\sum_{i=0}^\infty\binom{p}{i}J_{q+i-k+1}(a)u_2\\
&=\sum_{i=0}^\infty\binom{p}{i}a_{(q+i)}u_2
\end{split}
\end{equation}
and finally
\begin{equation}\label{eq-main-3}
\begin{split}
j_\infty(a\otimes z^p(z-1)^q(dz)^{1-k})u_3&=-\theta\Bigl(\sum_{i=0}^\infty(-1)^i\binom{q}{i}a\otimes
\xi_\infty^{p+q-i}(d\xi_\infty)^{1-k}\Bigr)u_3\\
&=-\sum_{i=0}^\infty(-1)^i\binom{q}{i}\theta(J_{p+q-i-k+1}(a))u_3
\end{split}
\end{equation}
for all $a\in V_k$ and $p,q\in\Z$.

By \eqref{eq-main-1}--\eqref{eq-main-3}, 
the definition of the functional $F$, Proposition \ref{prop-dual} and Proposition \ref{prop-eigen}, the left-hand side
of \eqref{eq-main-0} is equal to
\begin{equation}
\begin{split}
&\sum_{i=0}^\infty(-1)^{q-i}\binom{q}{i}\left\langle (u_2)^0_{(\alpha+q-i)}a_{(p+i)}u_1, u_3\right\rangle\\
&\quad\qquad\qquad+\sum_{i=0}^\infty\binom{p}{i}\left\langle(a_{(q+i)}u_2)_{(\alpha+p-i)}^0u_1,u_3\right\rangle\\
&\qquad\qquad\qquad\qquad\qquad-\sum_{i=0}^\infty(-1)^i\binom{q}{i}\left\langle a_{(p+q-i)}(u_2)_{(\alpha-i)}^0u_1, u_3
\right\rangle
\end{split}
\end{equation}
where $\alpha=|u_1|+|u_2|-|u_3|+k-2-p-q$,
which vanishes by $\eqref{eq-borcherds}$.
Hence \eqref{eq-main-0} is proved.

\medskip
\noindent
{\bf(Step 2)} We now prove that $I_x(-,z)\in I\fusion{M^2}{M^1}{D(M^3)}$.

\medskip
Since $M^i=\bigoplus_{n=0}^\infty M^i_{(h_i+n)}\in\mathscr{C}_{k_i}\,(i=1,2,3)$
we have
\begin{equation}\label{eq-main-4}
\begin{split}
&z^{-L_0}u_1\otimes z^{-L_0}u_2\otimes z^{L_0}u_3=
\sum_{n_1=0}^{k_1}\sum_{n_2=0}^{k_2}\sum_{n_3=0}^{k_3}
\frac{(-1)^{n_1+n_2}}{n_1!n_2!n_3!}\\
&\qquad\qquad\qquad\qquad\times(L_0-|u_1|)^{n_1}u_1\otimes(L_0-|u_2|)^{n_2}u_2\otimes
(L_0-|u_3|)^{n_3}u_3\\
&\qquad\qquad\qquad\qquad\qquad\qquad\qquad\qquad\times z^{|u_3|-|u_1|-|u_2|}(\log z)^{n_1+n_2+n_3}
\end{split}
\end{equation}
for homogeneous $u_1\in M^1$, $u_2\in M^2$ and $u_3\in M^3$.
Then the left-hand side of \eqref{eq-main-def-2} is an element in 
$\C[z,z^{-1}]z^{-h_1-h_2+h_3}[\log z]$.
Therefore $\left\langle I_x(u_2,z)u_1,-\right\rangle=\left\langle x, z^{-L_0}u_1\otimes z^{-L_0}u_2\otimes z^{L_0}-\right\rangle$
gives  an element of the space
\begin{equation}
\Hom_\C(M_3,\C)[[z,z^{-1}]]z^{-h_1-h_2+h_3}[\log z], 
\end{equation}
which shows 
$I_x(u_2,z)\in \Hom_\C(M_1, D(M_3))[[z,z^{-1}]]z^{-h_1-h_2+h_3}[\log z]$.
Therefore we can write 
\begin{equation}
I_x(u_2,z)=\sum_{n=0}^d\sum_{\alpha\in\Z+h_1+h_2-h_3}(u_2)_{(\alpha)}^nz^{-\alpha-1}(\log z)^n.
\end{equation}

For fixed $u_1\in M^1_{(h_1+\ell_1)}$ and $u_2\in M^2_{(h_2+\ell_2)}$ with nonnegative integers $\ell_1$
and $\ell_2$,
we have by \eqref{eq-main-4}
\begin{equation}\label{eq-main-10}
\langle I_x(u_2,z)u_1, u_3\rangle=\sum_{n=0}^{k_1+k_2+k_3}\sum_{\ell=0}^\infty c_\ell^n z^{h_3-h_1-h_2-\ell_1-\ell_2+\ell}
(\log z)^n
\end{equation}
where $c_\ell^n$ are complex numbers.
The \eqref{eq-main-10} implies that $(u_2)_{(\alpha)}^nu_1=0$ for $\alpha>h_3-h_1-h_2+\ell_1+\ell_2-1$.
Hence $I_x(-,z)$ satisfies the truncation condition.

In order to prove $L_{-1}$-derivative property,
we first note that
\begin{equation}\label{eq-main-11}
\begin{split}
&\left\langle x, j_0(\omega\otimes z(dz)^{-1})z^{-L_0}u_1\otimes z^{-L_0}u_2\otimes z^{L_0}u_3\right\rangle\\
&\quad\qquad+\left\langle x, z^{-L_0}u_1\otimes j_1(\omega\otimes z(dz)^{-1})z^{-L_0}u_2\otimes z^{L_0}u_3 \right\rangle\\
&\quad\qquad\qquad+\left\langle x, z^{-L_0}u_1\otimes z^{-L_0}u_2\otimes j_\infty(\omega\otimes z(dz)^{-1})z^{L_0}u_3 
\right\rangle=0.
\end{split}
\end{equation}
The left-hand side of \eqref{eq-main-11} turns to be
\begin{equation}\label{eq-deri}
\begin{split}
&\left\langle x,L_0z^{-L_0}u_1\otimes z^{-L_0}u_2\otimes z^{L_0}u_3\right\rangle\\
&\quad\qquad\qquad\qquad+\left\langle x, z^{-L_0}u_1\otimes (L_0+L_{-1})z^{-L_0}u_2\otimes z^{L_0}u_3 \right\rangle\\
&\quad\qquad\qquad\qquad\qquad\qquad\qquad-\left\langle x, z^{-L_0}u_1\otimes z^{-L_0}u_2\otimes L_0z^{L_0}u_3 \right\rangle.
\end{split}
\end{equation}
From now on each term in \eqref{eq-deri} is simplified.
Let us consider the second term of \eqref{eq-deri}.
Since $[L_{-1},L_0]=-L_{-1}$, we have
\begin{equation}
\left\langle x, z^{-L_0}u_1\otimes L_{-1}z^{-L_0}u_2\otimes z^{L_0}u_3 \right\rangle
=z\left\langle I_x(L_{-1}u_2,z)u_1, u_3\right\rangle,
\end{equation}
which shows
\begin{equation}\label{eq-deri1}
\begin{split}
z\left\langle I_x(u_2,z)u_1,u_3\right\rangle
=&-\left\langle x,L_0z^{-L_0}u_1\otimes z^{-L_0}u_2\otimes 
z^{L_0}u_3\right\rangle\\
&-\left\langle x, z^{-L_0}u_1\otimes L_0z^{-L_0}u_2\otimes z^{L_0}u_3 \right\rangle\\
&+\left\langle x, z^{-L_0}u_1\otimes z^{-L_0}u_2\otimes L_0z^{L_0}u_3 \right\rangle.
\end{split}
\end{equation}
The first term of \eqref{eq-deri1} can be calculated to be
\begin{equation}\label{eq-deri2}
\begin{split}
&\left\langle x,L_0z^{-L_0}u_1\otimes z^{-L_0}u_2\otimes z^{L_0}u_3\right\rangle\\
&=\sum_{n_1=0}^{k_1-1}\sum_{n_2=0}^{k_2}\sum_{n_3=0}^{k_3}
\frac{(-1)^{n_1+n_2}}{n_1!n_2!n_3!}
z^{-|u_1|-|u_2|+|u_3|}(\log z)^{n_1+n_2+n_3}\\
&\quad\times\left\langle x,(L_0-|u_1|)^{n_1+1}u_1\otimes (L_0-|u_2|)^{n_2}u_2\otimes (L_0-|u_3|)^{n_3}u_3\right\rangle\\
&\qquad\qquad+|u_1|\left\langle x, z^{-L_0}u_1\otimes z^{-L_0}u_2\otimes z^{L_0}u_3\right\rangle.
\end{split}
\end{equation}
Similarly, the second term of \eqref{eq-deri1} becomes
\begin{equation}\label{eq-deri3}
\begin{split}
&\left\langle x,z^{-L_0}u_1\otimes L_0z^{-L_0}u_2\otimes z^{L_0}u_3\right\rangle\\
&=\sum_{n_1=0}^{k_1}\sum_{n_2=0}^{k_2-1}\sum_{n_3=0}^{k_3}
\frac{(-1)^{n_1+n_2}}{n_1!n_2!n_3!}z^{-|u_1|-|u_2|+|u_3|}(\log z)^{n_1+n_2+n_3}\\
&\quad\times\left\langle x,(L_0-|u_1|)^{n_1}u_1\otimes (L_0-|u_2|)^{n_2+1}u_2\otimes (L_0-|u_3|)^{n_3}u_3\right\rangle\\
&\qquad\qquad+|u_2|\left\langle x, z^{-L_0}u_1\otimes z^{-L_0}u_2\otimes z^{L_0}u_3\right\rangle.
\end{split}
\end{equation}
Finally, the third term is
\begin{equation}\label{eq-deri4}
\begin{split}
&\left\langle x,z^{-L_0}u_1\otimes z^{-L_0}u_2\otimes L_0z^{L_0}u_3\right\rangle\\
&=\sum_{n_1=0}^{k_1}\sum_{n_2=0}^{k_2}\sum_{n_3=0}^{k_3-1}
\frac{(-1)^{n_1+n_2}}{n_1!n_2!n_3!}z^{-|u_1|-|u_2|+|u_3|}(\log z)^{n_1+n_2+n_3}\\
&\quad\times\langle x,(L_0-|u_1|)^{n_1}u_1\otimes (L_0-|u_2|)^{n_2}u_2\otimes (L_0-|u_3|)^{n_3+1}u_3\rangle\\
&\qquad\qquad-|u_3|\left\langle x, z^{-L_0}u_1\otimes z^{-L_0}u_2\otimes z^{L_0}u_3\right\rangle.
\end{split}
\end{equation}
In all of the calculations given above we have used the fact that 
each $M^i$ is an object in $\mathscr{C}_{k_i}$.
By using \eqref{eq-deri1}--\eqref{eq-deri4} we obtain
\begin{align}
&\frac{d}{dz}\left\langle x, z^{-L_0}u_1\otimes z^{-L_0}u_2\otimes z^{L_0}u_3\right\rangle\\
&\qquad =-z^{-1}\left\langle x,L_0z^{-L_0}u_1\otimes z^{-L_0}u_2\otimes z^{L_0}u_3\right\rangle\notag\\
&\qquad\qquad -z^{-1}\left\langle x, z^{-L_0}u_1\otimes (L_0)z^{-L_0}u_2\otimes z^{L_0}u_3\right\rangle\notag\\
&\qquad\qquad\qquad+z^{-1}\left\langle x, z^{-L_0}u_1\otimes z^{-L_0}u_2\otimes L_0z^{L_0}u_3\right\rangle,\notag
\end{align}
which shows
\begin{equation}
\frac{d}{dz}\left\langle x, z^{-L_0}u_1\otimes z^{-L_0}u_2\otimes z^{L_0}u_3\right\rangle
=\frac{d}{dz}\left\langle I_x(u_2,z)u_1,u_3\right\rangle.
\end{equation}
Hence we have proved the $L_{-1}$-derivative property 
\begin{equation}
I_x(L_{-1}u_2,z)u_1=\dfrac{d}{dz}I_x(u_2,z)u_1.
\end{equation}

Finally we will show \eqref{eq-borcherds}.
Since $x$ is a conformal block,
for any $p,\,q\in\Z$ and $a\in V_k$, we have

\begin{equation}
\begin{split}
&\langle x, j_0(a\otimes z^p(z-1)^q(dz)^{1-k})z^{-L_0}u_1\otimes z^{-L_0}u_2\otimes z^{L_0}u_3\rangle\\
&\quad+\langle x, z^{-L_0}u_1\otimes j_1(a\otimes z^p(z-1)^q(dz)^{1-k})z^{-L_0}u_2\otimes z^{L_0}u_3\rangle\\
&\qquad+\langle x, z^{-L_0}u_1\otimes z^{-L_0}u_2\otimes j_\infty(a\otimes z^p(z-1)^q(dz)^{1-k})z^{L_0}u_3\rangle=0.
\end{split}
\end{equation}
By \eqref{eq-main-1}--\eqref{eq-main-3}, \eqref{eq-involution} and $[L_0, a_{(n)}]=(k-n-1)a_{(n)}$ we have
\begin{equation}\label{eq-borcherds-prove}
\begin{split}
&\sum_{i=1}^\infty(-1)^{q+i}\binom{q}{i}z^{-p-i+k-1}\langle I_x(u_2,z)a_{(p+i)}u_1, u_3\rangle\\
&\quad+\sum_{i=1}^\infty\binom{p}{i}z^{-q-i+k-1}\langle I_x(a_{(q+i)}u_2,z)u_1, u_3\rangle\\
&\qquad-\sum_{i=0}^\infty(-1)^i\binom{q}{i}z^{-p-q+i+k-1}\langle a_{(p+q-i)}I_x(u_2, z)u_1, u_3\rangle
=0.
\end{split}
\end{equation}
Recall that  $I_x(u_2,z)u_1=\sum_{n=0}^d\sum_{\alpha\in\C}(u_2)_{(\alpha)}^nu_1z^{-\alpha-1}(\log z)^n$.
Then taking the coefficient
of $z^{-\alpha-p-q+k-2}(\log z)^n$ in \eqref{eq-borcherds-prove}  gives \eqref{eq-borcherds}.

\medskip
\noindent
{\bf (Step 3)} We will show that $F=x$ for any $x\in C^*(M_A,p_A)$ 
and that $I_{F}(-,z)=I(-,z)$ for $I(-,z)\in I\fusion{M^2}{M^1}{D(M^3)}$.

\medskip
Suppose that $I_x(u_2,z)=\sum_{n=0}^d\sum_{\alpha\in\C}(u_2)_{(\alpha)}^nz^{-\alpha-1}(\log z)^n$.
By \eqref{eq-main-4}, we have
\begin{equation}
\begin{split}
\langle F, u_1\otimes u_2\otimes u_3\rangle&=\langle I_x(u_2,1)u_1, u_3\rangle\\
&=\langle (u_2)_{(|u_1|+|u_2|-|u_3|-1)}^0u_1, u_3 \rangle\\
&=\langle x, u_1\otimes u_2\otimes u_3\rangle
\end{split}
\end{equation}
for any homogeneous $u_i\in M^i\,(i=1,2,3)$, which implies $F=x$.

Conversely, we see that
\begin{equation}\label{eq-425}
\begin{split}
&\langle I_{F}(u_2,z)u_1, u_3\rangle\\
&=\langle F, z^{-L_0}u_1\otimes z^{-L_0}u_2\otimes z^{L_0}u_3\rangle\\
&=\sum_{n_1=0}^{k_1}\sum_{n_2=0}^{k_2}\sum_{n_3=0}^{k_3}
\frac{(-1)^{n_1+n_2}}{n_1!n_2!n_3!}\\
&\quad\times\langle F, (L_0-|u_1|)^{n_1}u_1\otimes (L_0-|u_2|)^{n_2}u_2
\otimes (L_0-|u_3|)^{k_3}u_3\rangle\\
& \quad\times z^{-|u_1|-|u_2|+|u_3|}(\log z)^{n_1+n_2+n_3}\\
&=\sum_{n_1=0}^{k_1}\sum_{n_2=0}^{k_2}\sum_{n_3=0}^{k_3}
\frac{(-1)^{n_1+n_2}}{n_1!n_2!n_3!}\\
&\times\langle ((L_0-|u_2|)^{n_2}u_2)_{(\alpha)}^0(L_0-|u_1|)^{n_1}u_1, (L_0-|u_3|)^{n_3}u_3\rangle\\
&\quad\times z^{-|u_1|-|u_2|+|u_3|}(\log z)^{n_1+n_2+n_3}\\
&=\sum_{n_1=0}^{k_1}\sum_{n_2=0}^{k_2}\sum_{n_3=0}^{k_3}
\frac{(-1)^{n_1+n_2}}{n_1!n_2!n_3!}\\
&\times\langle (L_0-|u_3|)^{n_3} ((L_0-|u_2|)^{n_2}u_2)_{(\alpha)}^0(L_0-|u_1|)^{n_1}u_1, u_3\rangle\\
&\quad\times z^{-|u_1|-|u_2|+|u_3|}(\log z)^{n_1+n_2+n_3}
\end{split}
\end{equation}
where $\alpha=|u_1|+|u_2|-|u_3|-1$.
On the other hand, by Lemma \ref{lemma-fund}, we have
\begin{equation}\label{eq-426}
\begin{split}
&\sum_{n_1+n_2+n_3=k}
\frac{(-1)^{n_1+n_2}}{n_1!n_2!n_3!} (L_0-|u_3|)^{n_3} ((L_0-|u_2|)^{n_2}u_2)_{(\alpha)}^0(L_0-|u_1|)^{n_1}u_1\\
&=\sum_{n_1+n_2+n_3=k}
\frac{(-1)^{n_1+n_2}}{n_1!n_2!n_3!} \sum_{\ell=0}^{n_3}\binom{n_3}{\ell}\ell!(x_1+x_2)^{n_3-\ell}x_1^{n_1}x_2^{n_2}
(u_2)_{(\alpha)}^{\ell}u_1\\
&=\sum_{n_1+n_2+n_3=k}
\frac{(-1)^{n_1+n_2}}{n_1!n_2!n_3!} \sum_{\ell=0}^{n_3}\binom{n_3}{\ell}\ell!(x_1+x_2)^{n_3-\ell}x_1^{n_1}x_2^{n_2}
(u_2)_{(\alpha)}^{\ell}u_1\\
&=\sum_{\ell=0}^{k}\sum_{n_3=\ell}^{k}\Bigl(\sum_{n_1+n_2=k-n_3}\frac{(-1)^{n_1+n_2}}{n_1!n_2!(n_3-\ell)!}
(x_1+x_2)^{n_3-\ell}x_1^{n_1}x_2^{n_2}
(u_2)_{(\alpha)}^{\ell}u_1\Bigr)\\
&=\sum_{\ell=0}^{k}\sum_{n_3=0}^{k-\ell}\Bigl(\sum_{n_1+n_2=k-\ell-n_3}\frac{(-1)^{n_1+n_2}}{n_1!n_2!n_3!}
(x_1+x_2)^{n_3}x_1^{n_1}x_2^{n_2}
(u_2)_{(\alpha)}^{\ell}u_1\Bigr)\\
&=\sum_{\ell=0}^{k}\Bigl(\sum_{n_1+n_2+n_3=k-\ell}\frac{(-1)^{n_1+n_2}}{n_1!n_2!n_3!}
(x_1+x_2)^{n_3}x_1^{n_1}x_2^{n_2}
(u_2)_{(\alpha)}^{\ell}u_1\Bigr)\\
&=\sum_{\ell=0}^k \frac{1}{(k-\ell)!}(-x_1-x_2+x_1+x_2)^{k-\ell} (u_2)_{(\alpha)}^{\ell}u_1\\
&=(u_2)_{(\alpha)}^{k}u_1.
\end{split}
\end{equation}
Therefore, by combining \eqref{eq-425} and \eqref{eq-426} we obtain
\begin{equation}
\langle I_{F_I}(u_2,z)u_1, u_3\rangle=\langle I(u_2,z)u_1, u_3\rangle
\end{equation}
for homogeneous $u_i\in M^i\,(i=1,2,3)$.
The theorem is proved.

\end{document}